\documentclass[twoside,11pt]{amsart}

\newcommand{\ra}{\rightarrow}

\newcommand{\BC}{\mathbb{C}}
\newcommand{\BR}{\mathbb{R}}
\newcommand{\CP}{\mathbb{CP}}

\newcommand{\M}{\mathfrak{M}}

\newcommand{\TryPackage}[3]{\IfFileExists{#1.sty}{\usepackage{#1}#2}{#3}
}
\TryPackage{mathrsfs}{\renewcommand{\mathcal}{\mathscr}}{%
        \TryPackage{eucal}{}{}}

\newcommand{\ZZ}{{\mathbb Z}}
\newcommand{\RR}{{\mathbb R}}

\usepackage{latexsym}
\usepackage{amsmath}
\usepackage{amsthm}
\usepackage{amscd}
\usepackage{xypic} 
\usepackage{amsfonts}
\usepackage{amssymb}
\usepackage{graphicx}
\usepackage{psfrag}

\newtheorem{df}{Definition}
\newtheorem{thm}[df]{Theorem}
\newtheorem{cor}[df]{Corollary}
\newtheorem{lem}[df]{Lemma}
\newtheorem{prop}[df]{Proposition}

\newtheorem{conj}[df]{Conjecture}

\begin{document}

\title[Symplectic 4-manifolds near  the BMY line]{Symplectic 4-manifolds with arbitrary fundamental group near the Bogomolov-Miyaoka-Yau line}

\author{Scott Baldridge}
\author{Paul Kirk}

\thanks{The first   author  gratefully acknowledges support from the
NSF  grant DMS-0507857. The second  author  gratefully acknowledges
support from the NSF  grant DMS-0202148.}

\address{Department of Mathematics, Louisiana State University \newline
\hspace*{.375in} Baton Rouge, LA 70817}
\email{\rm{sbaldrid@math.lsu.edu}}

\address{Mathematics Department, Indiana University \newline
\hspace*{.375in} Bloomington, IN 47405}
\email{\rm{pkirk@indiana.edu}}

\maketitle

\begin{abstract} In this paper we construct a family of symplectic 4--manifolds with positive signature for any given fundamental group $G$ that approaches the BMY line.    The family is used to show that one cannot hope to do better than than the BMY inequality in finding a lower bound for the function $f=\chi+b\sigma$ on the class of all minimal symplectic 4-manifolds with a given fundamental group. \end{abstract}


\section{Introduction}

Let $\chi(S)$ and $\sigma(S)$ denote the Euler characteristic and signature of a 4-manifold respectively.  Minimal complex surfaces $S$ of general type satisfy $c_1^2(S)>0$, $\chi(S)>0$ and $$2\chi_h(S)-6\leq c_1^2(S) \leq 9\chi_h(S)$$ where $c_1^2(S)=2\chi(S)+3\sigma(S)$ and $\chi_h(S)=\frac14(\chi(S)+\sigma(S))$.
The second inequality is usually referred to as the Bogomolov-Miyaoka-Yau inequality.  Finding symplectic (or K\"{a}hler) 4-manifolds on or near the BMY line has a long and interesting history (c.f. \cite{MT}, \cite{Ch1}, \cite{Ch2}, \cite{PPX},  \cite{sw:4man_kirby_calc}, \cite{stip}, \cite{stip2}).  In particular this means looking for symplectic 4-manifolds with positive signature.

All known examples of symplectic 4--manifolds on the BMY line have large fundamental groups.  In fact, if $S$ is a complex surface differing from $\CP^2$ (the only known simply-connected example on the BMY line), the equality $c^2_1(S)=9\chi_h(S)$ holds if and only if the unit disk $D^4=\{(z_1,z_2)\in \BC^2 \ | \ |z_1|^2+ |z_2|^2\leq 1\}$ covers $S$  \cite{Hi1}, \cite{Y}, \cite{My}.  This implies  for such $S$ that $|\pi_1(S)|=\infty$.  The goal has generally been to produce examples that fill in the geography with respect to $(c_1^2,\chi_h)$.  In this paper we are interested what can be said for a given fundamental group.

Stipsicz  constructed simply connected symplectic 4--manifolds $C_n$ for which $c_1^2(C_n)/\chi_h(C_n) \ra 9$ as $n\ra \infty$ in \cite{stip2}.  Our main theorem generalizes this result to any fundamental group.

\begin{thm}
\label{thm:mainthm}
Let  $G$ have a   presentation with $g$
generators
$x_1,\cdots, x_g$ and $r$ relations $w_1,\cdots,w_r$.  For each integer $n>1$, there exists a  symplectic 4--manifold $M(G,n)$ with fundamental group $G$ with Euler characteristic
$$\chi(M(G,n))= 75n^2+256n+130+12(g+r+1),$$ and
signature
$$\sigma(M(G,n))= 25n^2-68n-78-8(g+r+1).$$
\end{thm}


Our interest in this question  developed while investigating pairs $(a,b)\in \BR^2$ for which the function $f=a\chi+b\sigma$ has a lower bound on the class of symplectic manifolds with a given fundamental group (see \cite{BK}). In that  article we  considered the following.

Fix a finitely presented group $G$ and let $\M$ denote either the class
$\M(G)$ of closed symplectic 4-manifolds with fundamental group $G$ or
the class $\M^{min}(G)$ of minimal, closed symplectic 4-manifolds with
fundamental group $G$.

 For $b\in \RR$, define $f_{\M}(b)\in \RR\cup \{-\infty\}$ to be the infimum
$$f_{\M}(b)=\inf_{M\in\M}\{ \chi(M)+b\sigma(M)\}.$$
(In \cite{BK} we considered the infimum $f_\M(a,b)$ of $a\chi+b\sigma$ on $\M$ and showed that
if $a\leq 0$ the infimum is $-\infty$. Thus we restrict to $f_\M(1,b)$, which  we more compactly denote by $f_\M(b)$ in the present article.)

We showed in \cite{BK} that the set
$$D_\M=\{b\ | f_\M(b)\ne -\infty\}$$ (the {\em domain of $f_\M$}) is an interval satisfying
$$[-1,1]\subset D_{\M(G)}\subset (-\infty,1]\text{ and }[-1,\tfrac32]\subset D_{\M^{\min}(G)}\subset (-\infty,\tfrac32].$$
The upper bounds are sharp; in fact $1\in D_{\M(G)}$  and $\frac32\in D_{\M^{min}(G)}$.

We are interested  in the value of the left endpoint  $e_G$  of $D_{\M(G)}$, which is an intriguing invariant of a group $G$. (It may or may not be contained in $D_{\M(G)}$.)   Since $e_G\leq -1$,  a  straightforward argument shows that $e_G$ is also  the left endpoint of  $D_{\M^{min}(G)}$.

 In \cite{BK} we observed that the  results of Stipcitz gives a better lower bound (than $-\infty$) when $G$ is the trivial group, and so a consequence of  the result of this article is an  extension to all  $G$.    In fact Theorem 1 easily implies the following corollary.
\begin{cor} For any  finitely presented  group $G$,
$$D_{\M(G)}\subset [-3,1]\text{ and }  D_{\M^{\min}(G)}\subset [-3,\tfrac32].$$
\qed \end{cor}

The BMY  inequality
$c_1^2\leq 9\chi_h$ is equivalent to   $f_{\M^{min}(G)}(-3)\ge 0$
provided $G$ is not a surface group.  Hence the  BMY conjecture and Corollary 2  together
imply that $e_G=3$. Thus, a weaker form of the BMY conjecture
could be stated as follows.

\begin{conj}[Weak BMY Conjecture]  For each finitely presented group $G$,  $e_G=-3$.
\end{conj}

\section{The construction}

We use the following notation. If $X$ and $Y$ are symplectic 4-manifolds containing symplectic surfaces $F_X\subset X$ and $F_Y\subset Y$ such that $F_X^2+F_Y^2=0$, then the symplectic sum   (\cite{symp:gompf:construct_symp_man}) of $X$ and $Y$ along $F_X$ and $F_Y$ will be  denoted by
$$X\#_{F_X,F_Y} Y.$$
Recall that topologically $X\#_{F_X,F_Y} Y$ is obtained by removing tubular neighborhoods of $F_X$ and $F_Y$ and identifying the resulting boundaries in a fiber-preserving way.  Moreover, if $E_X\subset X$ (resp. $E_Y\subset Y$) is a symplectic surface intersecting $F_X$ once transversally (resp. intersecting $F_Y$ transversally), then the fiber sum can be constructed so that (the connected  sum)
$E_X\#E_Y$ is a symplectic surface in $X\#_{F_X,F_Y}Y$.

\subsection{The first piece: symplectic manifolds with given fundamental group}

The following theorem was proven in  \cite{BK}.

\begin{thm} \label{thm:man_with_prescribed_fund_group} Let  $G$ have a   presentation with $g$
generators
$x_1,\cdots, x_g$ and $r$ relations $w_1,\cdots,w_r$. Then
   there exists a minimal symplectic 4-manifold $M(G)$ with $\pi_1(M(G))\cong G$,  Euler
characteristic
$\chi(M(G))=12(g+r+1)$, and signature $\sigma(M(G))=-8(g+r+1)$.
\end{thm}

We will use the following observation. The manifold $M(G)$ constructed in Theorem \ref{thm:man_with_prescribed_fund_group} is obtained by taking fiber sums of a certain base manifold with $g+r+1$ copies of the basic elliptic surface $E(1)$. Since $E(1)$ admits a singular fibration with symplectic generic fibers and 6 cusp fibers (which are simply connected), so does $E(1)-F$, where $F$ denotes the generic fiber in $E(1)$ along which the fiber sum giving $M$ is constructed.   Thus each $M(G)$ contains a symplectic torus $T_0$ such  that the induced homomorphism $\pi_1(T_0)\to \pi_1(M(G))$ is trivial.

\subsection{The second piece: symplectic manifolds near the BMY line}

In  \cite{stip2}, Stipsicz proved the following theorem.

\begin{prop}[Stipsicz] \label{prop:stip} For each non-negative integer $n$, there exists a symplectic 4--manifold $X(n)$ which admits a genus-(15n+1) Lefschetz fibration with a section $T_{n+2}$ of genus $(n+2)$ and self-intersection $-(n+1)$.  Furthermore, $X(n)$ can be equipped with a symplectic structure such that $T_{n+2}$ is a symplectic submanifold.  The projection map $X(n)\ra T_{n+2}$  induces an isomorphism on  fundamental groups.  The Euler characteristic of $X(n)$ is $\chi(X(n))=75n^2+180n+12$ and the signature is $\sigma(X(n))=25n^2-60n-8$.
\end{prop}

Denote by $F_{15n+1}\subset X(n)$ a fixed generic fiber of $X(n)$. This is a symplectic surface with trivial normal bundle.

\subsection{The third piece: a simply connected manifold}

Gompf constructs a symplectic 4-manifold $S_{1,1}$ in \cite[Lemma 5.5]{symp:gompf:construct_symp_man} which contains a disjoint pair $T,F$ of symplectically embedded surfaces $T$ of genus one and $F$ of genus two,  with trivial normal bundles such that $S_{1,1}-(T\cup F)$ is simply connected. Thus the symplectic sum $A$ of two copies $S_{1,1}$  along the genus two surfaces
$$A= S_{1,1}\#_{F,F}S_{1,1}$$
contains a pair of disjointly embedded symplectic tori $T_1\cup T_2\subset A$ with trivial normal bundles so that
the complement $A-(T_1\cup T_2)$ is simply connected. Since $S_{1,1}$ has Euler characteristic $23$ and signature $-15$, $\chi(A)=50$ and $\sigma(A)=-30$.

The  manifold $A$ has a useful property,  whose proof is a  simple application of the Seifert-Van Kampen  theorem.
\begin{prop} \label{pieceA}Suppose $B$ and $C$ are symplectic 4-manifolds containing symplectic  tori $i_B:T_B\subset B$   and  $i_C:T_C\subset C$   with trivial normal bundles.

Let  $D=B\#_{T_B,T_1} A\#_{T_2,T_C}C$ be the  fiber sum of $B,A$, and $C$. Then
$$\pi_1(D)=\big( \pi_1(B)/N((i_B)_*(\pi_1(T_B) )\big)\star\big( \pi_1(C)/N((i_C)_*(\pi_1(T_C)\big)$$
 where  $\star$ denotes free product and $N(H)$ denotes the normal closure of  a subgroup $H$.\qed
\end{prop}

\subsection{The fourth piece:  ``Elbows''}

Let $T$ be a torus and $\{a,b\}$ a pair of smoothly embedded loops forming a symplectic basis of $\pi_1T$.  Let $\varphi:T \ra T$  be the Dehn twist around $a$. The mapping torus $Y_\phi$ fibers over
$S^1$ with fiber $T$. Let $t_1:S^1\to Y_\phi$ denote a section.   Taking a product of $Y_\phi$ with $S^1$ yields a symplectic 4-manifold $Y_\phi\times S^1$ (this is just Thurston's manifold from \cite{thurston}) which fibers over a torus with symplectic torus fibers.   Moreover, the symplectic structure can be chosen so that the section $t_1\times \text{id}:S^1\times S^1\to Y_\phi\times S^1$ is symplectic. Denote by $s_1:S^1\to \{p\}\times S^1\subset Y_\phi\times S^1$ the loop representing the second factor.

Note that $Y_\phi\times S^1$ contains a torus $T'=b\times s_1$, where $b$ is the curve described above in the fiber of $Y_\phi$. The torus $T'$ is homologically non-trivial by the Kunneth theorem, since $b$ is non-trivial in $H_1(Y_\phi)$, and is Lagrangian with respect to the symplectic structure on $Y_\phi\times S^1$. Thus the symplectic structure on $Y_\phi\times S^1$ can be perturbed slightly to make $T'$ symplectic. Note moreover that $T'$ is disjoint from the section $t_1\times s_1:S^1\times S^1\to Y_\phi\times S^1$ since we can assume that $t_1$ intersects the fiber containing $b$ in a point which does not lie on $b$.  The tubular neighborhood of $T'$ in $Y_\phi\times S^1$ is trivial since $b$ can isotoped off itself in a fiber  of $Y_\phi\to S^1$.  Similarly the
tubular neighborhood of the section $t_1\times s_1$ is trivial since $t_1$ can be pushed off itself in $Y_\phi$.

Define $Elb(n)$ to be the fiber sum  $Elb(n) = (Y_\phi\times S^1) \#_{T,T^2} (T^2\times \Sigma_{n-1})$.  The fiber sum   can be carried out so that the sections of $ Y_\phi\times S^1\to  S^1\times S^1$ and $T^2\times \Sigma_{n-1}\to \Sigma_{n-1}$ yield a symplectic section of
the resulting fibration $Elb(n)\to\Sigma_n$.  Thus $Elb(n)$ contains a disjoint pair of symplectic surfaces with trivial normal bundles,
a torus $T'=b\times s_1$ and a genus $n$ surface, the image of the section,  which we denote by $ D_n$.

Letting $t_2,s_2,\cdots, t_n,s_n$ denote the generators of $\pi_1(\Sigma_{n-1})$, one computes
$$\begin{array}{lcl} \pi_1(Elb(n)) = \langle a,b, t_1,s_1, \dots, t_n, s_n  & | &
a \text{ central },[b,t_1]=a, \\
 && [b,t_i]=1 \text{ for } i>1,[b,s_i]=1 \text{ for all } i, \\
 &&  \prod_{i=1}^{n} [t_i,s_i]=1 \rangle.\end{array}$$

 The inclusion of $T'$ into  $Elb(n)$  takes the generators of $\pi_1T'$ to $b$ and  $s_1$, and the inclusion of $D_n$ takes the  standard basis to $t_1,s_1,\cdots,t_n,s_n$.  The Euler characteristic and signature of $Elb(n)$ both vanish.

The manifold $Elb(n)-D_n$ is a punctured torus fibration over $\Sigma_n$, and hence has a  presentation with  the same generators and  all the same relations  except  that one no longer has  $a$ commuting with  $b$. i.e. $a$ commutes with  all  generators  except  $b$.

\subsection{The fifth  piece: an elliptic surface}

 We find a symplectically embedded surface $J$ of genus $ n+3 $  and self-intersection $ n+1 $ in the elliptic surface $E(n+5)$  such that   $E(n+5)- J$ is simply connected as follows.  Consider $n+3$ copies of the generic fiber and one copy of the section in  a fibration $E(n+5)\ra \CP^1$ with $6(n+5)$ cusp fibers.  The section and fibers are symplectic with regards to the symplectic structure on the elliptic fibration $E(n+5)$.  Resolve the $n+3$ transverse double points (\cite{symp:gompf:construct_symp_man}) to get a symplectically embedded surface $J$ of genus $n+3$ and self-intersection $n+1$ (the fiber hits the section once and that section has self-intersection $-(n+5)$).
The complement $E(n+5)-J$ is simply connected because $E(n+5)$ has a simply connected fiber which intersects $J$ in one point: the normal circle of a tubular neighborhood of $J$ is nullhomotopic in $E(n+5)-J$.

\subsection{Putting the pieces together}

We begin by a modification of Stipcizs's construction.  Let $Z(n)$ be the fiber sum of
$Elb(15n+1)$ and $X(n)$ along $D_{15n+1}\subset Elb(15n+1)$ and the fiber $F_{15n+1}$ of the Lefschetz fibration $X(n)\to \Sigma_{n+2}$
$$Z(n)=Elb(15n+1)\#_{D_{15n+1},F_{15n+1}} X(n).$$
The fiber sum can  be constructed so that the fiber $T\subset Elb(15n+1)$ and the section $T_{n+2}\subset X(n)$ add to yield a symplectic surface of genus $n+3$,  $K_{n+3}=T\#T_{n+2}\subset Z(n)$ \cite{symp:gompf:construct_symp_man}.     The important property of  $Z(n)$ is that it contains a symplectic torus $T'$, since $D_{15n+1}$ and $T'$ are disjoint.

The fundamental group of $Z(n)$ is easily  computed, since  $Elb(15n+1)-D_{15n+1}$  is  a fiber bundle with punctured torus fibers and $X(n)-F_{15n+1}$ is a Lefschetz fibration over  a punctured genus $n+2$ surface with at least one simply connected fiber. Using the Seifert-Van Kampen  theorem and Novikov additivity one obtains the following.

\begin{lem} The fundamental  group of $Z(n)$   is the free product of $\ZZ$ with generator $b$, and a genus $n+2$ surface group generated by  $x_i,y_i$:
$$\pi_1(Z(n))=\ZZ b\star \langle  x_i,y_i, i=1,\cdots , n+2 |\prod[x_i,y_i]=1\rangle.$$

Moreover, $Z(n)$ contains a disjoint pair of symplectic  surfaces, $T'\cup K_{n+3}\subset Z(n)$
satisfying $[T']^2=0$, and
$[K_{n+3}]^2=-n-1$. The induced homomorphism $\pi_1(T')\to \pi_1 (Z(n)) $ is the  map
$$\langle a,s_1\ | \ [a,s_1]\rangle \to \pi_1 Z(n) \ a\mapsto  a, s_1\mapsto 1.$$
The induced homomorphism $\pi_1(K_{n+3})\to \pi_1(Z(n))$ is the  map
$$\langle a,b, x_1,y_1,\cdots, x_{n+2},y_{n+2}\ | \ [a,b]\prod[x_i,y_i]=1\rangle \to \pi_1(Z(n))
$$
$$a\mapsto  a, b\mapsto 1, x_i\mapsto x_i, y_i\mapsto y_i.$$

Moreover, $\chi(Z(n))=75n^2  +240n+12 $ and $\sigma(Z(n))=25n^2 -60n-8$. \qed
\end{lem}

The symplectic sum of $Z(n)$ with $E(n+5)$ along $J$,  $Z(n)\#_{K_{n+3},J}E(n+5)$
is a simply connected symplectic 4-manifold  containing  a torus $T_1$ with trivial normal  bundle and appropriate Euler characterstic and signature. We take symplectic sum of this manifold with $A$ to obtain an example with a torus whose complement is simply connected.

Define $W(n)$ to be the symplectic sum
$$W(n)=A\#_{T_1,T'}Z(n)\#_{K_{n+3},J}E(n+5).$$ Then since $\pi_1(A-(T_2\cup T_2))=1$, the following proposition follows straightforwardly.
\begin{prop} The symplectic manifold $W(n)$ is simply connected and contains  a symplectic torus $T_2\subset W(n)$ with trivial normal bundle so that $\pi_1(W(n)-T_2)=1$. It has Euler characterstic
 $\chi(W(n))=75n^2 +256 n +130$ and signature  $\sigma(W(n))=25n^2 -68n -78  .$  \qed
\end{prop}

We can now prove Theorem \ref{thm:mainthm}.

\begin{proof}[Proof of Theorem \ref{thm:mainthm}]
The symplectic sum
$$M(G,n)=M(G)\#_{T_0,T_2}W(n)$$
has fundamental group $G$  by Proposition \ref{pieceA}. The calculations of $\chi(M(G,n))$ and $\sigma(M(G,n))$ are routine.
\end{proof}







\end{document}